\documentclass[12pt,a4paper]{article}
\usepackage[T1]{fontenc}
\usepackage{amsmath,amsthm,mathtools}
\usepackage{newtxtext,newtxmath}
\usepackage[margin=1in]{geometry}
\usepackage{microtype}
\usepackage[hidelinks]{hyperref}
\hypersetup{pdftitle={Power-Law Lower Bounds for Edges in Odd Induced Subgraphs},pdfauthor={Xiaoyu Liu}}
\allowdisplaybreaks

\newtheorem{theorem}{Theorem}[section]
\newtheorem{lemma}[theorem]{Lemma}
\newtheorem{corollary}[theorem]{Corollary}
\newtheorem{proposition}[theorem]{Proposition}

\title{Power-Law Lower Bounds for Edges in Odd Induced Subgraphs}
\author{Xiaoyu Liu\thanks{Corresponding author. School of Mathematics and Statistics, Fuzhou University, Fuzhou 350108, China. Email: \texttt{2143084225@qq.com}.}}
\date{}

\begin{document}
\maketitle

\begin{abstract}
For a finite simple graph $G$, let $g(G)$ denote the maximum number of edges in an induced subgraph of $G$ in which every vertex has odd degree. Motivated by Scott's question of whether $g(G)\geq cm$ for every graph with $m$ edges, we prove that $g(G)\geq cm^{4/7}$ for every finite simple graph with $m\geq1$ edges, where $c>0$ is an absolute constant. For bipartite graphs, we obtain the stronger bound $g(B)\geq cm^{2/3}$.
\end{abstract}



\section{Introduction}

All graphs considered in this paper are finite and simple. A graph is called odd (even, respectively) if every vertex has odd (even, respectively) degree. Gallai's theorem states that the vertex set of every graph can be partitioned into two sets, each inducing an even subgraph. In studies of lower bounds on the order of odd induced subgraphs, the underlying graph is usually assumed to have no isolated vertices.

For a graph $G$, let $f_o(G)$ denote the maximum number of vertices in an odd induced subgraph of $G$. Let $g(G)$ denote the maximum number of edges in an odd induced subgraph of $G$. We allow the empty set in both definitions.
In particular, $g(G)=0$ when $G$ has no edges.

Lower bounds on the order of odd induced subgraphs have been studied extensively. For graphs without isolated vertices, Caro [2] established a lower bound of $c\sqrt{n}$, where $c>0$ is an absolute constant.
Scott [7] improved this bound to $n/\log n$, and Ferber and Krivelevich [3] later obtained $n/10000$, thereby resolving the corresponding linear-order
conjecture. Subsequently, Zeng [9] further improved the bound to $n/5$ using an odd-cut construction.
More recently, Guo, Gutin, Hao, Shi, Zhang and Zhou [4] extended the odd-cut approach to prescribed degree parities, obtaining a lower bound of $2n/11$. They also established rank-based bounds for bipartite graphs, including $f_o(B)\geq65n/256$.

These results concern the order of odd induced subgraphs rather than their number of edges. The present paper focuses on the edge-count problem and establishes lower bounds on the number of edges in odd induced subgraphs. Scott [8] raised the following problem.

\refstepcounter{problem}\noindent\textbf{Problem 1 (Scott [8]).} Does there exist an absolute constant $c>0$ such that every finite simple graph $G$ with $m$ edges satisfies $g(G)\geq cm$?\par

The known linear lower bounds on the order of odd induced subgraphs yield only a square-root lower bound on the number of edges. Indeed, after deleting isolated vertices, let $n=|V(G)|$. Since every vertex in an odd induced subgraph has degree at least one, an odd induced subgraph of order $\Omega(n)$ contains at least $\Omega(n)$ edges. On the other hand, $m\leq\binom{n}{2}$ for a simple graph. Consequently, $g(G)=\Omega(\sqrt{m})$.

In this paper, we improve the exponent in the general lower bound from $1/2$ to $4/7$. We also establish lower bounds in terms of the minimum degree and for several special classes of graphs. Our main results are stated below.

\begin{theorem}
There exists an absolute constant $c>0$ such that every finite simple graph $G$ with $m\geq1$ edges satisfies $g(G)\geq cm^{4/7}$.
\end{theorem}

\begin{theorem}
There exists an absolute constant $c>0$ such that every finite simple graph $G$ with minimum degree $d\geq1$ satisfies $g(G)\geq cd^{4/3}$.
\end{theorem}

\begin{theorem}
There exists an absolute constant $c>0$ such that every bipartite graph $B$ with $m\geq1$ edges satisfies $g(B)\geq cm^{2/3}$.
\end{theorem}

We also obtain a lower bound for graphs whose neighborhoods have bounded chromatic number.

\begin{theorem}
	There exists an absolute constant $c>0$ with the following property. Let $G$ be a graph with $m$ edges such that $\chi(G[N(v)])\leq r$ for every $v\in V(G)$, where $r\geq1$. Then
	\[
	g(G)\geq \frac{c m^{3/5}}{r^{1/5}(\log(2+m))^{2/5}}.
	\]
\end{theorem}

In particular, Theorem 1.4 applies to triangle-free graphs and $K_{2,t}$-free graphs for $t\geq2$.

The proof begins with a fractional odd-subgraph lemma, obtained by combining Zeng's odd-cut constructions [9] with a domination result of Bollobás and Cockayne [1]. This lemma guarantees a probability distribution on the odd induced subgraphs of every graph without isolated vertices, under which each vertex is selected with probability at least $1/6$. Together with the weighted Caro--Wei inequality, it yields an independent-set extraction lemma for local cuts.

We next establish a parity-completion method based on a common pivot. By combining random selection with a binary linear system, we obtain odd induced subgraphs containing sufficiently many edges. A dyadic decomposition of the degree ranges then gives the minimum-degree bound. For bipartite graphs, complementary degree-layer estimates yield the stronger edge-count bound.

The paper is organized as follows. Section 2 establishes the fractional odd-subgraph lemma and the local-cut extraction lemma. Section 3 develops the parity-completion tools and proves Theorem 1.3. Section 4 is devoted to the proofs of Theorems 1.2 and 1.1, while Section 5 establishes Theorem 1.4 and its consequences for graphs with colorable neighborhoods. Section 6 discusses open problems.

For disjoint sets $A,B\subseteq V(G)$, let $E_G(A,B)$ denote the set of edges between $A$ and $B$, and write $e_G(A,B)=|E_G(A,B)|$. For $S\subseteq V(G)$, write $d_S(x)=|N_G(x)\cap S|$. For $F\subseteq E(G)$, let $d_F(x)$ denote the number of edges in $F$ incident with $x$, and put $e_F(A,B)=|F\cap E_G(A,B)|$. Subscripts are omitted when no ambiguity arises. All logarithms are natural. Throughout the paper, $c,C,c_1,C_1,\ldots$ denote positive absolute constants, which may vary from one occurrence to another.

\section{A fractional odd-subgraph lemma and local cuts}

We recall Gallai's theorem [5], which states that every graph admits a partition of its vertex set into two sets, each inducing an even subgraph.

\begin{lemma}[Zeng~{[9, Theorem 1.2]}]
	Let $H$ be a graph and let $V(H)=A\mathbin{\dot\cup}B$. Suppose that every vertex has an odd number of neighbors in the other part. Then $V(H)$ can be partitioned into two sets, each inducing an odd subgraph of $H$.
\end{lemma}

For a dominating set $D$ of a graph $G$, an \emph{external private neighbor} of $u\in D$ is a vertex $x\in V(G)\setminus D$ such that $N(x)\cap D=\{u\}$.

\begin{lemma}[Bollob\'as and Cockayne~\cite{BC}]
	Every graph without isolated vertices has a minimum dominating set in which every vertex has an external private neighbor.
\end{lemma}

Combining the two odd-cut constructions of Zeng~\cite[Section~3]{Zeng}, we obtain a probability distribution on the odd induced subgraphs of every graph without isolated vertices, such that each vertex is selected with probability at least $1/6$.

\begin{theorem}
	Let $H$ be a graph without isolated vertices. There exists a probability distribution on its odd induced subgraphs such that each vertex is selected with probability at least $1/6$.
	
	Consequently, for every function $w:V(H)\to\mathbb{R}_{\geq 0}$, there exists an odd induced subgraph $J$ satisfying $w(V(J))\geq\frac16w(V(H))$, where $w(S)=\sum_{x\in S}w(x)$ for $S\subseteq V(H)$.
\end{theorem}

\begin{proof}
Let $D$ be a minimum dominating set satisfying Lemma 2.2. For each $u\in D$, define $Q_u=\{x\in V(H)\setminus D:N(x)\cap D=\{u\}\}$. Each $Q_u$ is nonempty, and the sets $Q_u$ are pairwise disjoint. Set $Q=\bigcup_{u\in D}Q_u$ and $R=V(H)\setminus(D\cup Q)$.

Following Zeng~\cite[Section~3]{Zeng}, we use two odd-cut constructions to define probability distributions on odd induced subgraphs.

For the first distribution, choose one vertex $q_u$ uniformly from each $Q_u$. Consider the induced subgraph $H[D\cup\{q_u:u\in D\}]$, with the partition $D\mathbin{\dot\cup}\{q_u:u\in D\}$. Every $u\in D$ has exactly one neighbor across this partition, namely $q_u$, and every selected vertex $q_u$ has exactly one neighbor across the partition, namely $u$. Thus this is an odd cut. Apply the odd-cut lifting lemma and choose one of the two resulting odd induced subgraphs with probability $1/2$. Writing $t_u=|Q_u|$, a vertex $x\in D$ is present in the selected odd subgraph with probability $1/2$. For $x\in Q_u$, it must first be chosen as $q_u$, which has probability $1/t_u$, and then lie in the selected part, which has probability $1/2$. Hence its inclusion probability is $1/(2t_u)$. Vertices in $R$ are not used in this construction.

For the second distribution, choose a random set $Z\subseteq D$ by including each vertex independently with probability $1/2$, and put $R_Z=\{x\in R:d_Z(x)\text{ is odd}\}$. For each $u\in Z$, choose a largest set $T_u\subseteq Q_u$ such that $|T_u|+d_{R_Z}(u)\equiv1\pmod 2$. Thus either all of $Q_u$ is retained or exactly one vertex is deleted. If a deletion is required, choose the deleted vertex uniformly and independently for each $u$. In particular, each $x\in Q_u$ is retained with probability at least $(t_u-1)/t_u$ conditional on $u\in Z$.

Consider the induced subgraph $H[Z\cup R_Z\cup\bigcup_{u\in Z}T_u]$, with the partition $Z\mathbin{\dot\cup}(R_Z\cup\bigcup_{u\in Z}T_u)$. A vertex $u\in Z$ has cross-degree $d_{R_Z}(u)+|T_u|$, which is odd by the choice of $T_u$. A vertex $x\in R_Z$ has cross-degree $d_Z(x)$, which is odd by the definition of $R_Z$. Finally, every vertex of $T_u$ has exactly one neighbor in $Z$, namely $u$, because $T_u\subseteq Q_u$. Thus this is also an odd cut.

Apply the odd-cut lifting lemma and choose one of the two resulting odd induced subgraphs uniformly. A vertex $x\in D$ is selected with probability at least $1/4$: it must belong to $Z$, and then it is in the selected part with probability $1/2$. 

For $x\in R$, the set $N(x)\cap D$ is nonempty, since $D$ dominates $H$. The parity of $d_Z(x)$ is therefore odd with probability $1/2$; conditional on $x\in R_Z$, the vertex is selected with probability $1/2$. Its inclusion probability is at least $1/4$. For $x\in Q_u$, it must be that $u\in Z$, that $x$ is retained in $T_u$, and that it lies in the selected odd part. Its inclusion probability is consequently at least $(1/2)((t_u-1)/t_u)(1/2)=(t_u-1)/(4t_u)$.

Now use the first distribution with probability $1/3$ and the second with probability $2/3$. A vertex of $D$ is selected with probability at least $1/3$, and a vertex of $R$ with probability at least $1/6$. For $x\in Q_u$, its inclusion probability is at least $\frac13\cdot\frac{1}{2t_u}+\frac23\cdot\frac{t_u-1}{4t_u}=\frac16$. This proves the first assertion. Taking expectations of the selected weight gives the second assertion.
\end{proof}

\begin{corollary}
Every $n$-vertex graph without isolated vertices has an odd induced subgraph on at least $n/6$ vertices and containing at least $n/12$ edges.
\end{corollary}

\begin{proof}
Apply the preceding theorem with $w(x)=1$ for every vertex. There is an odd induced subgraph $J$ with $|V(J)|\geq n/6$. Since every vertex of $J$ has degree at least one, the handshaking lemma gives $2e(J)=\sum_{x\in V(J)}d_J(x)\geq|V(J)|\geq n/6$. Hence $e(J)\geq n/12$.
\end{proof}

We shall combine the fractional odd-subgraph lemma with the weighted Caro--Wei inequality [6].

\begin{lemma}
Let $H$ be a graph and let $w:V(H)\to\mathbb R_{\geq0}$. Then $H$ contains an independent set $I$ such that $w(I)\geq\sum_{x\in V(H)}\frac{w(x)}{d_H(x)+1}$. If $W=w(V(H))>0$, then $I$ may be chosen so that $w(I)\geq\frac{W^2}{W+\sum_{x\in V(H)}w(x)d_H(x)}$.
\end{lemma}

\begin{proof}
Take a uniformly random ordering of $V(H)$ and put $x$ into $I$ if it appears before every vertex of $N_H(x)$. The resulting set is independent, and $\Pr(x\in I)=1/(d_H(x)+1)$. Taking expectations gives an independent set whose weight is at least $\sum_x w(x)/(d_H(x)+1)$.

For the second assertion, the Cauchy--Schwarz inequality gives $W^2\leq(\sum_x w(x)/(d_H(x)+1))(\sum_x w(x)(d_H(x)+1))$. The second factor equals $W+\sum_xw(x)d_H(x)$. The first assertion therefore gives an independent set of weight at least the stated quotient.
\end{proof}

The following local-cut lemma is the basic extraction tool. Keeping the maximum designated degree is useful in the applications below.

\begin{lemma}
Let $U\subseteq V(G)$ and let $F\subseteq E_G(U,V(G)\setminus U)$, where $q=|F|\geq1$. Put $b=\max_{x\in V(G)\setminus U}d_F(x)$. For every $\tau>0$, either $g(G)\geq\tau$, or there is an independent set $I\subseteq V(G)\setminus U$ such that $e_F(U,I)\geq\frac{q^2}{72(q+b\tau)}$.
\end{lemma}

\begin{proof}
Let $V_F$ be the set of endpoints of the edges in $F$, and let $H_0=G[V_F]$. Every vertex of $H_0$ is incident with an edge of $F$, so $H_0$ has no isolated vertices. Define $w(x)=d_F(x)$ for $x\notin U$ and $w(x)=0$ for $x\in U$. Each edge of $F$ has exactly one endpoint outside $U$, and hence $w(V(H_0))=q$.

By the fractional odd-subgraph lemma, $H_0$ has an odd induced subgraph $J$ such that $W:=w(V(J))\geq q/6$. If $e(J)\geq\tau$, then $J$ is also an odd induced subgraph of $G$, and the first alternative holds. Otherwise, the weighted Caro--Wei inequality gives an independent set $I_0\subseteq V(J)$ such that $w(I_0)\geq\frac{W^2}{W+\sum_{x\in V(J)}w(x)d_J(x)}$.

Since $w(x)\leq b$ for every $x$ and $e(J)<\tau$, we have $\sum_{x\in V(J)}w(x)d_J(x)\leq b\sum_{x\in V(J)}d_J(x)=2b e(J)<2b\tau$. Also $q/6\leq W\leq q$.

Therefore, $w(I_0)\geq\frac{q^2}{36(q+2b\tau)}\geq\frac{q^2}{72(q+b\tau)}$.

Put $I=I_0\setminus U$. The set $I$ is independent and lies outside $U$. Removing the zero-weight vertices in $U$ does not change its weight. Since each edge of $F$ has exactly one endpoint outside $U$, every edge of $F$ incident with a vertex of $I$ is counted exactly once in $w(I)$. Thus $w(I)=e_F(U,I)$, as required.
\end{proof}

\section{Parity completion and bipartite graphs}

All matrix operations in this section are over $\mathbb F_2$.

Binary linear systems have also been used recently in prescribed-parity vertex-order problems; see Guo et al. [4]. Here we use a weighted kernel-selection argument in which the coordinates are weighted by their contributions to the cross-edge count. Thus the rank loss is measured against edge weight, which is essential for the edge-count objective considered in this paper.

We first establish a weighted selection lemma for the kernel of a binary matrix.

\begin{lemma}
Let $M$ be a binary matrix, and assign a weight $\omega_j\in[0,a]$ to each column $j$. There is a vector $z\in\ker M$ such that $\sum_j\omega_jz_j\geq\frac12(\sum_j\omega_j-a\operatorname{rank}M)$.
\end{lemma}

\begin{proof}
Choose a vector uniformly from $\ker M$. A coordinate $z_j$ is fixed on $\ker M$ if and only if the standard basis vector $e_j$ belongs to the row space of $M$. Indeed, a coordinate is fixed precisely when the corresponding coordinate functional vanishes on $\ker M$, which is equivalent to $e_j\in(\ker M)^\perp=\operatorname{row}(M)$.

Since the standard basis vectors are linearly independent, at most $\operatorname{rank}M$ coordinates can be fixed. Every fixed coordinate is zero, because $0\in\ker M$. For every nonfixed coordinate, there is a vector $u\in\ker M$ with $u_j=1$. Translation by $u$ pairs each kernel vector with another whose $j$-th coordinate is different. Thus a nonfixed coordinate equals $1$ with probability exactly $1/2$ under the uniform distribution on the kernel.

It follows that $\mathbb E[\sum_j\omega_jz_j]\geq\frac12(\sum_j\omega_j-a\operatorname{rank}M)$. Some vector in the kernel therefore attains at least this value.
\end{proof}

\begin{lemma}
Let $B$ be bipartite with bipartition $(P,Q)$. Let $z\in P$, and let $R\subseteq N_B(z)\cap Q$ have odd cardinality. Define $Q(R)=\{x\in P\setminus\{z\}:d_R(x)\text{ is odd}\}$. Then $g(B)\geq\frac12(e_B(R,Q(R))-|R|^2)$.
\end{lemma}

\begin{proof}
Let $M$ be the binary adjacency matrix with rows indexed by $R$ and columns indexed by $Q(R)$. Assign to column $x$ the weight $d_R(x)$. Every column weight is at most $|R|$, and $\operatorname{rank}M\leq|R|$. By the preceding kernel-selection lemma, there is a vector $\mathbf u\in\ker M$ whose support $T\subseteq Q(R)$ satisfies $e_B(R,T)\geq\frac12(e_B(R,Q(R))-|R|^2)$, because the sum of the column weights is exactly $e_B(R,Q(R))$.

The equation $M\mathbf u=0$ means that $d_T(r)$ is even for every $r\in R$. Moreover, every $x\in T$ belongs to $Q(R)$, so $d_R(x)$ is odd. Consider the induced subgraph $B[R\cup T\cup\{z\}]$. Each $r\in R$ has degree $1+d_T(r)$, each $x\in T$ has degree $d_R(x)$, and $z$ has degree $|R|$. All these degrees are odd. Thus this induced subgraph is odd, and the required inequality follows.
\end{proof}

We use the elementary parity identity $\Pr(\operatorname{Bin}(q,p)\text{ is even})=(1+(1-2p)^q)/2\geq1/2$ for every integer $q\geq0$ and $0\leq p\leq1/2$.

The next lemma gives a lower bound when two independent sets have a common pivot with the specified neighborhood.

\begin{lemma}
There exists an absolute constant $c>0$ with the following property. Let $X,Y$ be disjoint independent sets in a graph $G$, and let $z\notin X\cup Y$ satisfy $N_G(z)\cap(X\cup Y)=X$. Put $s=|X|\geq1$, $M=e_G(X,Y)$, and $\lambda=M/s$. Then $g(G)\geq c\min\{\lambda^2,M\}$.
\end{lemma}

\begin{proof}
There is nothing to prove if $M=0$. If $0<\lambda\leq16$, a single edge gives the result after decreasing $c$, since $\min\{\lambda^2,M\}\leq256$. Suppose that $\lambda>16$. Put $r=\min\{\lambda,s\}$ and $p=r/(64s)$, and choose a random set $R\subseteq X$ by selecting each vertex independently with probability $p$. Let $Q(R)=\{y\in Y:d_R(y)\text{ is odd}\}$, and let $\mathcal A$ be the event that $|R|$ is odd.

Fix an edge $xy\in E_G(X,Y)$ and condition on $x\in R$. The remaining vertices of $X$ split into the disjoint sets $A=(N_G(y)\cap X)\setminus\{x\}$ and $B=X\setminus N_G(y)$. Under this conditioning, $y\in Q(R)$ holds precisely when the number of selected vertices in $A$ is even. The event $\mathcal A$ requires the total number selected in $A\cup B$ to be even. Thus the two events occur together precisely when the selected numbers in both $A$ and $B$ are even. These counts are independent, and the parity identity above shows that each is even with probability at least $1/2$. Consequently, $\Pr(y\in Q(R),\mathcal A\mid x\in R)\geq1/4$.

Summing the contribution of each edge gives $\mathbb E[e_G(R,Q(R))\mathbf1_{\mathcal A}]\geq pM/4=r\lambda/256$. Also, $\mathbb E[|R|^2\mathbf1_{\mathcal A}]\leq\mathbb E|R|^2=sp(1-p)+(sp)^2\leq r/64+r^2/4096$. Since $r\leq\lambda$ and $\lambda>16$, we have $r/64+r^2/4096\leq5r\lambda/4096$. Subtracting the latter estimate from the former yields $\mathbb E[(e_G(R,Q(R))-|R|^2)\mathbf1_{\mathcal A}]\geq11r\lambda/4096$.

Hence there is an outcome for which $|R|$ is odd and $e_G(R,Q(R))-|R|^2\geq11r\lambda/4096$. Apply the parity-completion lemma to the induced bipartite graph $G[X\cup Y\cup\{z\}]$, with bipartition $(Y\cup\{z\},X)$ and pivot $z$. Since $r\lambda=\min\{\lambda^2,M\}$, the result follows.
\end{proof}

The next lemma is the main estimate for a single dyadic degree layer.

\begin{lemma}
There exists an absolute constant $c>0$ such that the following holds. Let $B$ be bipartite with bipartition $(X,Y)$, where $s=|X|\geq1$, and suppose that $k\leq d_X(y)<2k$ for every $y\in Y$, for some integer $k\geq1$. Put $M=e_B(X,Y)$ and $\lambda=M/s$. Then $g(B)\geq ck\lambda=ckM/s$.
\end{lemma}

\begin{proof}
If $M=0$, the assertion is immediate. We first prove the weaker estimate $g(B)\geq c_0k\min\{k,\lambda\}$. If $1\leq k\leq16$, then $k\min\{k,\lambda\}\leq k^2\leq256$, so the existence of an edge proves this estimate after decreasing $c_0$. Thus assume $k\geq17$.

Choose $x_0\in X$ with $t=d_B(x_0)\geq\lambda$, and put $W=N_B(x_0)$. Since every $y\in W$ has $k\leq d_X(y)<2k$, we have $e_B(W,X\setminus\{x_0\})=\sum_{y\in W}(d_X(y)-1)\geq(k-1)t$.

Let $r=\min\{k,t\}$, $p=r/(64t)$, and select each vertex of $W$ independently with probability $p$, obtaining a random set $R$. Define $Q(R)=\{x\in X\setminus\{x_0\}:d_R(x)\text{ is odd}\}$, and let $\mathcal A$ be the event that $|R|$ is odd.

For an edge $xy$ with $y\in W$ and $x\ne x_0$, condition on $y\in R$. The requirements $x\in Q(R)$ and $\mathcal A$ hold simultaneously precisely when the selected numbers in $N_B(x)\cap W\setminus\{y\}$ and $W\setminus N_B(x)$ are both even. These two sets are disjoint and partition $W\setminus\{y\}$, so the corresponding counts are independent. By the parity identity, the conditional probability is at least $1/4$. Summing over these edges gives $\mathbb E[e_B(R,Q(R))\mathbf1_{\mathcal A}]\geq p(k-1)t/4=r(k-1)/256$. Also, $\mathbb E[|R|^2\mathbf1_{\mathcal A}]\leq\mathbb E|R|^2\leq r/64+r^2/4096$.

Since $r\leq k$ and $k\geq17$, the difference of these two bounds satisfies
\[
\frac{r(k-1)}{256}-\frac r{64}-\frac{r^2}{4096}
\geq\frac{rk}{500}.
\]
Therefore, some outcome with $|R|$ odd satisfies $e_B(R,Q(R))-|R|^2\geq rk/500$. Applying the parity-completion lemma with pivot $x_0$ and bipartition $(X,Y)$, we obtain $g(B)\geq c_1rk$. Since $t\geq\lambda$, we have $r\geq\min\{k,\lambda\}$. This proves the weaker estimate.

If $\lambda\leq128k$, the weaker estimate implies $g(B)\geq ck\lambda$: this is immediate when $\lambda\leq k$, while for $k<\lambda\leq128k$ we use $k^2\geq k\lambda/128$. It remains to consider $\lambda>128k$.

When $k=1$, every vertex of $Y$ has exactly one neighbor in $X$. Since $M=s\lambda$, some vertex of $X$ has at least $\lambda$ neighbors in $Y$. Taking the largest odd number of these neighbors gives an induced star with at least $c\lambda$ edges, as required.

Now assume $k\geq2$. For $y\in Y$, put $E_y=\sum_{x\in N_B(y)}(d_B(x)-1)$. By double counting and the Cauchy--Schwarz inequality,
\[
\begin{aligned}
\sum_{y\in Y}E_y
&=\sum_{x\in X}d_B(x)(d_B(x)-1)&=\sum_{x\in X}d_B(x)^2-M
&\geq\frac{M^2}{s}-M&=M(\lambda-1).
\end{aligned}
\]
Since $|Y|\leq M/k$, there exists $y_0\in Y$ such that $E_{y_0}\geq k(\lambda-1)\geq k\lambda/2$.

Let $W=N_B(y_0)$, so $|W|<2k$. Select each vertex of $W$ independently with probability $p=1/64$, obtaining a random set $R$. Define $Q(R)=\{y\in Y\setminus\{y_0\}:d_R(y)\text{ is odd}\}$, and let $\mathcal A$ be the event that $|R|$ is odd.

For each edge $xy$ with $x\in W$ and $y\ne y_0$, condition on $x\in R$. The two requirements $y\in Q(R)$ and $\mathcal A$ hold simultaneously precisely when the selected numbers in $N_B(y)\cap W\setminus\{x\}$ and $W\setminus N_B(y)$ are both even. These sets are disjoint and partition $W\setminus\{x\}$, so the conditional probability is at least $1/4$. Summing over these edges gives $\mathbb E[e_B(R,Q(R))\mathbf1_{\mathcal A}]\geq pE_{y_0}/4\geq k\lambda/512$.

Moreover, since $|W|<2k$ and $p=1/64$, we have $\mathbb E[|R|^2\mathbf1_{\mathcal A}]\leq\mathbb E|R|^2\leq2pk+4p^2k^2=k/32+k^2/1024\leq17k^2/1024<k\lambda/1024$, where the last inequality follows from $\lambda>128k$. Consequently, some outcome with $|R|$ odd satisfies $e_B(R,Q(R))-|R|^2\geq k\lambda/1024$. Applying the parity-completion lemma with pivot $y_0$ and bipartition $(Y,X)$ gives $g(B)\geq ck\lambda$, as required.
\end{proof}

We now prove the bipartite theorem.

\begin{proof}[Proof of Theorem 1.3]
Delete isolated vertices and let $(X,Y)$ be a bipartition of $B$, where $a=|X|\leq|Y|$. Write $m=e(B)$ and $t=g(B)$. Partition $Y$ into the dyadic classes $Y_j=\{y\in Y:2^j\leq d_X(y)<2^{j+1}\}$, where $j$ ranges over nonnegative integers. Put $k_j=2^j$ and $m_j=e_B(X,Y_j)$.

By the degree-layer lemma, $t\geq c k_jm_j/a$, and hence $m_j\leq Ct a/k_j$. On the other hand, $m_j<2k_j|Y_j|$, so $|Y_j|>m_j/(2k_j)$. After deleting isolated vertices from $B[X\cup Y_j]$, at least $|Y_j|$ vertices remain. Corollary 2.4 therefore gives $t\geq c|Y_j|\geq cm_j/k_j$, and hence $m_j\leq Ct k_j$.

Combining the two estimates, we obtain $m_j\leq Ct\min\{k_j,a/k_j\}$. Since the $k_j$ are powers of two, the two geometric sums on either side of $\sqrt a$ are each at most $2\sqrt a$. Thus $\sum_j\min\{k_j,a/k_j\}\leq4\sqrt a$, and consequently $m=\sum_jm_j\leq Ct\sqrt a$.

Applying Corollary 2.4 to $B$ itself gives $t\geq cn$, where $n=|V(B)|$. Since $n\geq2a$, we have $t\geq ca$, and thus $a\leq Ct$. It follows that $m\leq Ct^{3/2}$, which yields $t\geq cm^{2/3}$. This proves Theorem 1.3.
\end{proof}

\section{The general edge bound}

We next prove the minimum-degree theorem.

\begin{proof}[Proof of Theorem 1.2]
Fix the constants in the local-cut lemma and the two common-pivot lemmas. Choose an absolute constant $\eta>0$ sufficiently small so that the two contradictions below hold. For all sufficiently large $d$, put $\tau=\eta d^{4/3}$. We may assume that $\tau\geq d$ and $d\geq8$. Suppose, for a contradiction, that $g(G)<\tau$.

Choose a vertex $v$ of degree $d$, put $N=N_G(v)$, and let $X$ be a maximum independent set in $G[N]$. Write $s=|X|$. Apply Lemma 2.6 to the star with center $v$ and leaf set $N$, taking $U=\{v\}$ and threshold $\tau$. Here $q=d$ and $b=1$. Since $g(G)<\tau$, the second alternative gives an independent set $I\subseteq N$ with $|I|=e_F(\{v\},I)\geq d^2/[72(d+\tau)]$. By the maximality of $X$ and the assumption $\tau\geq d$, we obtain $s\geq |I|\geq d^2/[72(d+\tau)]\geq c d^2/\tau$.

Let $O=V(G)\setminus(N\cup\{v\})$. We claim that $e_G(X,O)<sd/4$. Suppose otherwise. Choose $F\subseteq E_G(X,O)$ with $q=\lfloor sd/4\rfloor$ edges. Since $d\geq8$ and $s\geq1$, we have $sd/8\leq q\leq sd/4$. Moreover, every vertex outside $X$ has designated degree at most $s$, so the parameter $b$ in Lemma 2.6 satisfies $b\leq s$.

Apply Lemma 2.6 with $U=X$ and threshold $\tau$. Since $g(G)<\tau$, there exists an independent set $Y\subseteq O$ such that $M:=e_G(X,Y)\geq q^2/[72(q+b\tau)]$. Using $q\geq sd/8$, $q\leq sd/4$, $b\leq s$, and $\tau\geq d$, we have $q+b\tau\leq sd/4+s\tau\leq2s\tau$. Therefore $M\geq c s d^2/\tau$.

Writing $\lambda=M/s$, we have both $s\geq cd^2/\tau$ and $\lambda\geq cd^2/\tau$. The vertex $v$ is adjacent to every vertex of $X$ and to no vertex of $Y$. Since $X$ and $Y$ are independent, the common-pivot lemma applies with pivot $v$. It follows that $g(G)\geq c\lambda\min\{\lambda,s\}\geq cd^4/\tau^2=c\eta^{-3}\tau$. Choosing $\eta$ sufficiently small contradicts $g(G)<\tau$. This proves the claim.

Since $X$ is independent and every vertex of $X$ is adjacent to $v$, each vertex of $X$ has exactly $d-1$ other neighbors, all lying in $N\setminus X$ or $O$. Thus $e_G(X,N\setminus X)=s(d-1)-e_G(X,O)>3sd/4-s\geq sd/2$, where the last inequality uses $d\geq8$.

Apply Lemma 2.6 again, now with $U=X$, $F=E_G(X,N\setminus X)$, and threshold $\tau$. We have $sd/2\leq q\leq sd$ and $b\leq s$. Since $\tau\geq d$, $q+b\tau\leq sd+s\tau\leq2s\tau$. Consequently, Lemma 2.6 gives an independent set $Y\subseteq N\setminus X$ such that $M:=e_G(X,Y)\geq c s d^2/\tau$.

We now partition $Y$ according to the degrees of its vertices into $X$. For each nonnegative integer $j$, let $Y_j=\{y\in Y:2^j\leq d_X(y)<2^{j+1}\}$, and put $k_j=2^j$, $m_j=e_G(X,Y_j)$, and $\lambda_j=m_j/s$. Since $Y_j$ is independent in $G[N]$ and $X$ is a maximum independent set in $G[N]$, we have $|Y_j|\leq s$. Hence $\lambda_j<2k_j$.

On the other hand, applying the degree-layer lemma to $G[X\cup Y_j]$ gives $g(G)\geq c k_j\lambda_j$. The assumption $g(G)<\tau$ therefore implies $\lambda_j\leq C\tau/k_j$. Combining the two estimates, we obtain $M/s=\sum_j\lambda_j\leq C\sum_j\min\{k_j,\tau/k_j\}$.

Since the $k_j$ are powers of two, splitting at $\sqrt\tau$ and summing the two geometric tails gives $\sum_j\min\{k_j,\tau/k_j\}\leq4\sqrt\tau$. Hence $M/s\leq C\sqrt\tau$.

On the other hand, the preceding lower bound on $M$ gives $M/s\geq cd^2/\tau$. Consequently, $d^2/\tau\leq C\sqrt\tau$, which implies $\tau\geq c d^{4/3}$. This contradicts $\tau=\eta d^{4/3}$ when $\eta$ is sufficiently small.

The argument covers all sufficiently large $d$. For bounded positive $d$, an edge gives $g(G)\geq1$, and decreasing the constant completes the proof.
\end{proof}

\begin{proof}[Proof of Theorem 1.1]
Delete all isolated vertices from $G$, and put $n=|V(G)|$. Since $m\geq1$, we have $n\geq2$. By Corollary 2.4, $g(G)\geq n/12$.

Set $h=m/n$. Starting with $G$, repeatedly delete a vertex whose degree in the current graph is less than $h$, as long as such a vertex exists. This process must leave a nonempty induced subgraph $H$. Indeed, if every vertex were deleted, then each original edge would be counted exactly once, when its first endpoint was deleted. Since every deleted vertex has current degree less than $h$, summing the degrees at the deletion steps would give $m<nh=m$, a contradiction.

Consequently, $H$ is nonempty and satisfies $\delta(H)\geq h$. Since $h>0$ and degrees are integers, we also have $\delta(H)\geq1$. By Theorem 1.2 and the fact that $H$ is an induced subgraph of $G$, we obtain $g(G)\geq g(H)\geq c(m/n)^{4/3}$.

Combining this estimate with $g(G)\geq n/12$, and taking the weighted geometric mean with weights $4/7$ and $3/7$, respectively, gives
\[
\begin{aligned}
g(G)
&\geq \left(\frac n{12}\right)^{4/7}
       \left(c\left(\frac mn\right)^{4/3}\right)^{3/7}\\
&\geq c'm^{4/7}.
\end{aligned}
\]
This proves Theorem 1.1.
\end{proof}

\section{Graphs with colorable neighborhoods}

The logarithmic factor in the following result arises from selecting one degree layer.

\begin{proposition}
There exists an absolute constant $c>0$ such that the following holds. Let $G$ have minimum degree at least $d\geq1$, and suppose that $\chi(G[N(v)])\leq r$ for every vertex $v$. Then $g(G)\geq c d^{3/2}/(\sqrt r\,\log(2+d))$.
\end{proposition}

\begin{proof}
It is enough to consider sufficiently large $d$, since for bounded positive $d$, an edge gives the result after decreasing $c$. Let $D=\delta(G)\geq d$, and choose a vertex $v$ of degree $D$. Since $G[N(v)]$ is $r$-colorable, one of its color classes is an independent set $U$ of size $s:=|U|\geq D/r$.

Let $Z=V(G)\setminus(U\cup\{v\})$. Since $U$ is independent and every vertex of $U$ is adjacent to $v$, each vertex of $U$ has at least $D-1$ neighbors in $Z$. Hence $e_G(U,Z)\geq s(D-1)\geq sD/2$.

For each nonnegative integer $j$, put $k_j=2^j$ and $Z_j=\{z\in Z:k_j\leq d_U(z)<2k_j\}$. Since $d_U(z)\leq s\leq D$, there are at most $C\log(2+D)$ nonempty classes. Let $J$ be the number of nonempty classes. The classes partition all edges between $U$ and $Z$, so $e_G(U,Z)=\sum_j e_G(U,Z_j)$. Since $e_G(U,Z)\geq sD/2$ and $J\leq C\log(2+D)$, at least one class, with lower degree $k$, contains $q:=e_G(U,Z_j)\geq csD/\log(2+D)$ edges from $U$.

Put $t=g(G)$ and let $F=E_G(U,Z_j)$ be the set of these $q$ edges. Apply Lemma 2.6 with designated set $U$ and threshold $t+1$. The alternative $g(G)\geq t+1$ is impossible by the definition of $t$. Every vertex in $Z_j$ has designated degree less than $2k$, so the parameter $b$ in Lemma 2.6 is at most $2k$. 

Hence $q+b(t+1)\leq q+2k(t+1)\leq2(q+k(t+1))$. Lemma 2.6 therefore gives an independent set $Y\subseteq Z_j$ such that, writing $M=e_G(U,Y)$, we have $M\geq c q^2/[q+k(t+1)]$. The set $Y$ lies in $Z_j$ because all endpoints of $F$ outside $U$ lie there.

The graph $G[U\cup Y]$ is bipartite, and every vertex of $Y$ has degree between $k$ and $2k-1$ toward $U$. Applying Lemma 3.4 to this induced subgraph, and using $g(G)\geq g(G[U\cup Y])$, gives $t\geq ckM/s$. Moreover, $|Y|\geq M/(2k)$, since every vertex of $Y$ has fewer than $2k$ neighbors in $U$. After deleting isolated vertices from $G[U\cup Y]$, at least $|Y|$ vertices remain. Corollary 2.4 thus gives $t\geq c|Y|\geq cM/k$.

Combining these two estimates, we obtain $t^2\geq c(kM/s)(M/k)=cM^2/s$, and hence $t\geq cM/\sqrt{s}$.

We distinguish two cases. If $q\geq k(t+1)$, then $q+k(t+1)\leq2q$, and the preceding lower bound on $M$ gives $M\geq cq$. Therefore, $t\geq cq/\sqrt{s}\geq cD\sqrt{s}/\log(2+D)\geq cD^{3/2}/(\sqrt r\,\log(2+D))$. Here we used $q\geq csD/\log(2+D)$ and $s\geq D/r$.

If $q<k(t+1)$, then $q+k(t+1)<2k(t+1)$, so the same estimate for $M$ gives $M\geq cq^2/[k(t+1)]$. Combining this with $t\geq ckM/s$, we obtain $t\geq cq^2/[s(t+1)]$. Since $t\geq1$, we have $t+1\leq2t$, and hence $t^2\geq cq^2/s$. Thus $t\geq cq/\sqrt{s}$, which yields the same lower bound as in the first case.

A direct differentiation shows that $x^{3/2}/\log(2+x)$ is increasing for $x\geq1$. Since $D\geq d$, we conclude that $g(G)\geq c d^{3/2}/(\sqrt r\,\log(2+d))$. This proves the proposition.
\end{proof}

\begin{proof}[Proof of Theorem 1.4]
Delete all isolated vertices from $G$, and put $n=|V(G)|$. By Corollary 2.4, $g(G)\geq n/12$.

If $m/n<1$, then $n>m$, and hence $g(G)>m/12$. Since $r\geq1$ and $m^{3/5}/(\log(2+m))^{2/5}\leq m$ for $m\geq1$, the required bound follows after decreasing the constant.

We may therefore assume that $m/n\geq1$. As in the proof of Theorem 1.1, repeatedly deleting vertices of current degree less than $m/n$ leaves a nonempty induced subgraph $H$ with $\delta(H)\geq m/n$. The neighborhood-coloring hypothesis is inherited by induced subgraphs, since $N_H(v)\subseteq N_G(v)$ for every $v\in V(H)$. Put $D=\delta(H)$. Proposition 5.1 gives $g(G)\geq g(H)\geq cD^{3/2}/(\sqrt r\,\log(2+D))$.

Since $D\geq m/n\geq1$ and the function $x^{3/2}/\log(2+x)$ is increasing for $x\geq1$, we have $D^{3/2}/\log(2+D)\geq(m/n)^{3/2}/\log(2+m/n)\geq(m/n)^{3/2}/\log(2+m)$. Thus $g(G)\geq c r^{-1/2}(m/n)^{3/2}/\log(2+m)$.

Taking the weighted geometric mean of this bound with weight $2/5$ and $g(G)\geq n/12$ with weight $3/5$, the powers of $n$ cancel and yield
\[
\begin{aligned}
g(G)
&\geq \left(c r^{-1/2}
       \frac{(m/n)^{3/2}}{\log(2+m)}\right)^{2/5}
       \left(\frac n{12}\right)^{3/5}\\
&\geq c' r^{-1/5}
       \frac{m^{3/5}}{(\log(2+m))^{2/5}}.
\end{aligned}
\]
This proves Theorem 1.4.
\end{proof}

\begin{corollary}
There exists an absolute constant $c>0$ such that every triangle-free graph $G$ with $m\geq1$ edges satisfies $g(G)\geq cm^{3/5}/(\log(2+m))^{2/5}$. For every fixed $t\geq2$, the same conclusion, with a constant depending on $t$, holds for $K_{2,t}$-free graphs.
\end{corollary}

\begin{proof}
The neighborhood of every vertex in a triangle-free graph is independent. If $G$ is $K_{2,t}$-free, every neighborhood has maximum degree at most $t-1$, and hence is $t$-colorable. The result follows from Theorem 1.4.
\end{proof}

\section{Open Problems}

Scott's linear edge-count question, stated as Problem 1, remains open. Theorem 1.1 gives the lower bound $g(G)\geq cm^{4/7}$ for general graphs, while Theorem 1.2 gives a lower bound of order $d^{4/3}$ in terms of the minimum degree. The exponent $4/3$ in Theorem 1.2 comes from the balance $d^2/\tau\leq C\sqrt\tau$: the first quantity is the normalized edge count supplied by the local-cut extraction, while the second bounds the total contribution of the dyadic degree layers. Thus, removing logarithmic losses alone does not improve the exponent. A stronger general bound may require a parity-completion argument that retains internal edges on one side of the cut, rather than first replacing both sides by independent sets.

The results proved here lead to several natural questions. Can the exponent $4/7$ in Theorem 1.1 be improved? More specifically, does there exist an absolute constant $\varepsilon>0$ such that every finite simple graph with $m$ edges satisfies $g(G)\geq c m^{4/7+\varepsilon}$ for some absolute constant $c>0$? It is also natural to ask whether the minimum-degree bound in Theorem 1.2 can be improved beyond the exponent $4/3$.

For bipartite graphs, Theorem 1.3 gives an exponent of $2/3$. It remains to determine whether this exponent is best possible, or whether a stronger lower bound can be obtained. Finally, for graphs whose neighborhoods are $r$-colorable, can the logarithmic factor in Theorem 1.4 be reduced or removed while retaining a comparable dependence on $r$?

These questions indicate that the edge-count problem for odd induced subgraphs remains open in several natural settings.

\section*{Declaration on the Use of AI}
The author used ChatGPT (OpenAI) to assist in discussing proof strategies, checking proofs, and improving the exposition.


\begin{thebibliography}{99}
\bibitem{BC} B. Bollob\'as and E. J. Cockayne, Graph-theoretic parameters concerning domination, independence, and irredundance, \emph{J. Graph Theory} \textbf{3} (1979), 241--249.

\bibitem{Caro} Y. Caro, On induced subgraphs with odd degrees, \emph{Discrete Math.} \textbf{132} (1994), 23--28.

\bibitem{FK} A. Ferber and M. Krivelevich, Every graph contains a linearly sized induced subgraph with all degrees odd, \emph{Adv. Math.} \textbf{406} (2022), 108534.

\bibitem{Guo} Q. Guo, G. Gutin, Y. Hao, Y. Shi, Y. Zhang and Y. Zhou, Bounds on odd and odd-even induced subgraphs, arXiv:2608.02339v2 [math.CO], 2026.

\bibitem{Lovasz} L. Lov\'asz, \emph{Combinatorial Problems and Exercises}, 2nd ed., North-Holland Publishing Co., Amsterdam, 1993.

\bibitem{STY} S. Sakai, M. Togasaki and K. Yamazaki, A note on greedy algorithms for the maximum weighted independent set problem, \emph{Discrete Appl. Math.} \textbf{126} (2003), 313--322.

\bibitem{Scott92} A. D. Scott, Large induced subgraphs with all degrees odd, \emph{Combin. Probab. Comput.} \textbf{1} (1992), 335--349.

\bibitem{Scott01} A. D. Scott, On induced subgraphs with all degrees odd, \emph{Graphs Combin.} \textbf{17} (2001), 539--553.

\bibitem{Zeng} Q. Zeng, Large odd induced subgraphs via odd cuts, arXiv:2607.23266v1 [math.CO], 2026.
\end{thebibliography}
\end{document}